# Some Probabilistic and Statistical Properties of a Random Coefficient Autoregressive Model

A. Bouchemella and A. Bibi

Badji Mokhtar - Annaba University

abdelhalimgbs@yahoo.fr

November 2, 2018

**Abstract**

A statistical inference for random coefficient first-order autoregressive model $[RCAR(1)]$ was investigated by P.M. ROBINSON (1978) in which the coefficients varying over individuals. In this paper we attempt to generalize this result to random coefficient autoregressive model of order $p$ $[RCAR(p)]$. The stationarity condition will derived for this model.

## 1  Introduction

We can envisage several variants of the random coefficients autoregressive models. Most these variants are constructed while raising the different hypotheses progressively on the structure of the coefficients. One is interested then here in the $RCAR(p)$ model whose coefficients vary by individual, while trying to study the statistical inference of this model one being based on Robinson's article (1978).

We consider the $RCAR(p)$ model, satisfying

$$Y_t(\omega) = \sum_{k=1}^{p} A_k(\omega) Y_{t-k}(\omega) + \varepsilon_t(\omega). \tag{1.1}$$

The model (1.1) becomes

$$\begin{bmatrix} Y_{t-p+1}(\omega) \\ Y_{t-p+2}(\omega) \\ \vdots \\ Y_t(\omega) \end{bmatrix} = \begin{bmatrix} 0 & 1 & 0 & \cdots & 0 \\ 0 & 0 & \ddots & \ddots & \vdots \\ 0 & 0 & \ddots & \ddots & 0 \\ \vdots & \vdots & \cdots & 0 & 1 \\ A_p(\omega) & A_{p-1}(\omega) & \cdots & \cdots & A_1(\omega) \end{bmatrix} \begin{bmatrix} Y_{t-p}(\omega) \\ Y_{t-p+1}(\omega) \\ \vdots \\ Y_{t-1}(\omega) \end{bmatrix} + \begin{bmatrix} 0 \\ 0 \\ \vdots \\ \varepsilon_t(\omega) \end{bmatrix} \tag{1.2}$$

or
$$\underline{Y}_t(\omega) = A(\omega)\underline{Y}_{t-1}(\omega) + \underline{\varepsilon}_t(\omega) \tag{1.3}$$

Where $\underline{Y}_t(\omega)$ is the observable random variable, the $t \in \{0,1,2,...,T\}$ and $\omega \in \{1,2,...,N\}$ indexes designate the time and the population of individuals respectively.

We suppose that :

(i) $\{\varepsilon_t(\omega)\}$ is an independent sequence of univariate random variables with mean zero and variance $\sigma_\omega^2$.

(ii) $A_i(\omega)$, $1 \leq i \leq p$ are independent random variables with $E\{A_1(\omega),...,A_p(\omega)\} = (\alpha_1,...,\alpha_p)'$.

(iii) $A_i(\omega)$, $1 \leq i \leq p$ are mutually independent of $\varepsilon_t(\omega)$ and $Y_t(\omega)$ for each $t$.

(iv) $\{A(\omega), \omega = 1,...,N\}$ is an independent sequence of $p \times p$ matrices with

$$E\{A(\omega)\} = A = \begin{bmatrix} 0 & 1 & 0 & \cdots & 0 \\ 0 & 0 & \ddots & \ddots & \vdots \\ 0 & 0 & \ddots & \ddots & 0 \\ \vdots & \vdots & \cdots & 0 & 1 \\ \alpha_p & \alpha_{p-1} & \cdots & \cdots & \alpha_1 \end{bmatrix}$$

(v) There is no non-zero $p \times 1$ constant vector $Z$ such that $Z'Y_t(\omega)$ is determined exactly as a linear function of $\{Y_{t-1}(\omega),...,Y_{t-p}(\omega)\}$.

(vi) $Y_t(\omega)$ is independent of $\varepsilon_s(\omega)$ for each $s > t$.

In the following section we interested to study the conditions of stationarity as well as the conditions of existence and uniqueness of the stationary solutions. In the section 2 while basing on Robinson's survey that made on a $RCAR(1)$ model, we try to have a generalization of the results gotten for the $RCAR(p)$ model.

## 2 Stationarity

Nicholls and Quinn [6] use some models which also belong to the class of autoregressive processes of which the coefficients varying with time, and they got the conditions of stationarity for these models.

Now we will find these conditions when the coefficients vary by groups of individuals under certain assumptions.

We note that $\mathcal{F}_t(\omega)$ is $\sigma$-field generated by a sequence $\{\varepsilon_s(\omega) \text{ and } A_k(\omega), 1 \leq k \leq p, s \leq t\}$ for each $\omega \in \{1,...,N\}$.

**Theorem 2.1** *The process $\{Y_t(\omega)\}$ $t = 0, 1, ...$ given by (1.1) from $t = 1$ is stationary if and only if*

$$\begin{cases} \beta_0 = \beta_1 \\ \\ \Gamma_\omega(0) = \Gamma_\omega(1) \end{cases}$$

Where $\beta_i = E\{\underline{Y}_i(\omega)\}$ and $\Gamma_\omega(i-j) = E\{\underline{Y}_i(\omega)\underline{Y}_j(\omega)\}$.

**Proof.** See [6]. ∎

**Theorem 2.2** *The process $\{Y_t(\omega)\}$; $t = 0, 1, ...$ given by (1.1) is stationary if and only if $\beta = E\{\underline{Y}_0(\omega)\}$ satisfies $A\beta = \beta$ and $\Gamma_\omega(0) = E\{\underline{Y}_0(\omega)\underline{Y}_0(\omega)\}$ satisfies the equation*

$$Vec[\Gamma_\omega(0)] = (A \otimes A) Vec[\Gamma_\omega(0)] + \Omega_\omega \tag{2.1}$$

*where $\otimes$ is the Kronecker product.*

**Lemma 2.1** *If the condition*

$$z^p - \alpha_1 z^{p-1} - \cdots - \alpha_p \neq 0 \text{ for } |z| \geq 1 \tag{2.2}$$

*is verified then the equation (2.1) has the unique positive definite matrix solution. This solution is explicitly given by*

$$Vec[\Gamma_\omega(0)] = \sum_{k=0}^{\infty} A^k \Omega_\omega A'^k \tag{2.3}$$

**Corollary 2.1** *For the process $\{Y_t(\omega)\}$ to be stationary it is necessary that a condition (2.2) is verified.*

**Remark 2.1** *The condition (2.2) is equivalent that the eigenvalues of the matrix $A$ are less than unity in modulus.*

Anděl [1] has shown that under certain conditions on $\varepsilon_t(\omega)$ and $A_i(\omega)$, it is possible to get stronger properties for $\mathcal{F}_t(\omega)$-measurable solution $\{Y_t(\omega)\}$ to (1.1), in which casd random variables.

**Theorem 2.3** *In order that there exist a stationary $\mathcal{F}_t(\omega)$-measurable solution $\{Y_t(\omega)\}$ to (1.1) satisfying assumptions (i)-(iv), it is necessary that $\sum_{j=0}^{r}(A \otimes A)^j Vec[\Omega_\omega]$ converge as $r \to \infty$.*

*When $(A \otimes A)$ does not have a unit eigenvalue, this latter condition is both necessary and sufficient, and there is a unique stationary solution $\{Y_t(\omega)\}$ obtained from*

$$\underline{Y}_t(\omega) = \eta_t(\omega) + \sum_{j=1}^{\infty} A^j(\omega) \eta_{t-j}(\omega) \tag{2.4}$$

*Where*

$$\eta_t(\omega) = \begin{bmatrix} 0 & 0 & \cdots & 0 \\ \vdots & \ddots & \ddots & \vdots \\ 0 & \cdots & 0 & 0 \\ 0 & \cdots & 0 & \varepsilon_t(\omega) \end{bmatrix}$$

*Moreover when $\{\varepsilon_t(\omega)\}$ and $\{A_i(\omega)\}_{i=1}^{p}$ are also identically distributed sequences. Then this solotion is also strictly stationary and ergodic.*

**Proof.** See [6]. ∎

# 3 Statistical inference for RCAR(p)

In this section we will retain the same assumptions $(i)$-$(vi)$ set forth above, adding that

$(vii)$ $\{\varepsilon_t(\omega)\}$ and $\{A_i(\omega)\}_{i=1}^{p}$ are also identically distributed sequences.

So under the conditions $(i)$-$(vii)$ and in accordance with the theorem of existence and uniqueness above there is the unique solution measurable strictly stationary and ergodic.

**Remark 3.1** *In what follows one agrees that :*

$\Gamma_\omega(u) = E\{\underline{Y}_t(\omega)\underline{Y}_{t-u}(\omega)'/\mathcal{F}_t(\omega)\}$ *a.s i.e. indicate the conditional covariance matrix.*

$\Upsilon_\omega(u) = E\{\underline{Y}_t(\omega)\underline{Y}_{t-u}(\omega)'\}$ *indicate the uncondinal covariance matrix.*

Then we have

$$E\{\underline{Y}_t(\omega)/\mathcal{F}_t(\omega)\} = 0 \quad a.s \tag{3.1}$$

$$\Gamma_\omega(u) = E\{\underline{Y}_t(\omega)\underline{Y}_{t-u}(\omega)'/\mathcal{F}_t(\omega)\} \quad a.s \tag{3.2}$$

where $\Gamma_\omega(.)$ is $p \times p$ the covariance matrix.

For $u = 0$ we have

$$\begin{aligned}
\Gamma_\omega(0) &= E\{\underline{Y}_t(\omega)\underline{Y}'(\omega)/\mathcal{F}_t(\omega)\} \\
&= E\{[A(\omega)\underline{Y}_{t-1}(\omega) + \underline{\varepsilon}_t(\omega)][A(\omega)\underline{Y}_{t-1}(\omega) + \underline{\varepsilon}_t(\omega)]'/\mathcal{F}_t(\omega)\} \\
&= A(\omega)E\{\underline{Y}_{t-1}(\omega)\underline{Y}'_{t-1}(\omega)/\mathcal{F}_t(\omega)\}A'(\omega) + E\{\underline{\varepsilon}_t(\omega)\underline{\varepsilon}'_t(\omega)/\mathcal{F}_t(\omega)\}
\end{aligned}$$

then

$$\Gamma_\omega(0) = A(\omega)\Gamma_\omega(0)A'(\omega) + \Omega_\omega \quad a.s$$

The solution of this equation is

$$vec[\Gamma_\omega(0)] = [I_{p^2} - A^{\otimes 2}(\omega)]^{-1} vec[\Omega_\omega] \quad a.s \tag{3.3}$$

$$\Gamma_\omega(u) = A(\omega)\Gamma_\omega(u-1), \quad u = 1, 2, ... \; a.s$$

by recurrence we obtain

$$\Gamma_\omega(u) = A^u(\omega)\Gamma_\omega(0), \quad u = 1, 2, ... \; a.s \tag{3.4}$$

$$Vec[\Gamma_\omega(u)] = (I_p \otimes A^u(\omega))Vec[\Gamma_\omega(0)] \quad a.s$$

$$Vec[\Gamma_\omega(u)] = I_p \otimes A^u(\omega)[I_{p^2} - A^{\otimes 2}(\omega)]^{-1}Vec[\Omega_\omega] \quad a.s \tag{3.5}$$

Since the matrix $\Gamma_\omega(0)$ is non singular, it may be identified as

$$\begin{aligned}
A(\omega) &= E\{\underline{Y}_t(\omega)\underline{Y}'_{t-1}(\omega)/\mathcal{F}_t(\omega)\}\left[E\{\underline{Y}_t(\omega)\underline{Y}'_t(\omega)/\mathcal{F}_t(\omega)\}\right]^{-1} \\
&= \Gamma_\omega(1)\Gamma_\omega^{-1}(0) \tag{3.6}
\end{aligned}$$

Consider $\hat{A}_T(\omega)$ given by

$$\hat{A}_T(\omega) = \left[\sum_{t=1}^{T} \underline{Y}_t(\omega) \underline{Y}'_{t-1}(\omega)\right] \left[\sum_{t=1}^{T} \underline{Y}_t(\omega) \underline{Y}'_t(\omega)\right]^{-1} \tag{3.7}$$

Then, under slightly stronger conditions above we have that, for $u \geq 0$

$$p \lim_{T \to \infty} \left( \left\| T^{-1} \sum_{t=1}^{T+u} \underline{Y}_t(\omega) \underline{Y}'_{t-u}(\omega) - E\{\underline{Y}_t(\omega) \underline{Y}'_{t-u}(\omega) / \mathcal{F}_t(\omega)\} \right\| \right) = 0 \tag{3.8}$$

for any fixed $u$ so that

$$p \lim_{T \to \infty} \left( \left\| \hat{A}_T(\omega) - A(\omega) \right\| \right) = 0 \tag{3.9}$$

We assume that $A(\omega)$ has all its eigenvalues within the unit circle. Then in the usual terminology $\{\underline{Y}_t(\omega)\}$ is an $AR(p)$ and $\{\underline{Y}_t(\omega)\}$ is ergodic in the sense that, if $\Upsilon_\omega(u) = E\{\underline{Y}_t(\omega) \underline{Y}'_{t-u}(\omega)\}$ then (3.8) and (3.9) become

$$p \lim_{T \to \infty} \left( \left\| T^{-1} \sum_{t=1}^{T+u} \underline{Y}_t(\omega) \underline{Y}'_{t-u}(\omega) - \Upsilon_\omega(u) \right\| \right) = 0, \ \forall u$$

$$p \lim_{T \to \infty} \left( \left\| \hat{A}_T(\omega) - A \right\| \right) = 0$$

where $A$ is a such matrix, which its eigenvalues are in the unit circle unit.

In the next subsection we will show that the moments of $A(\omega)$ can be identified in terms of $\Upsilon_\omega(u)$, then we will estimate the moments of $A(\omega)$ while studying their asymptotic behavior as $N \to \infty$ while $T$ stays fixed, under some conditions on the eigenvalues of the matrix $A(\omega)$.

But an alternative approach more evident believes $A(\omega)$ by one of the usual estimators for example estimator ordinary least squares, then it is estimated moments from, but this approach has several disadvantages such as:

Especially when $T$ is not large enough we can not estimate $A(\omega)$ and our approach seems more quicker and direct. Second the estimator of $A(\omega)$ will not have in fact eigenvalues in the unit circle unit. Third, the asymptotic properties of estimators of $\mu(u)$ obtained from the estimator of $A(\omega)$ are less easy and accessible than those obtained by our estimators.

## 3.1 Identifying moments of $\mathbf{A}(\omega)$

**Theorem 3.1** *If the matrix $A(\omega)$ has all its eigenvalues within the unit circle, then :*

$$vec[\Upsilon_\omega(0)] = \sum_{v=0}^{\infty} \underline{\mu}_{2v} Vec[\Omega_\omega] \tag{3.10}$$

$$Vec[\Upsilon_\omega(u)] = \sum_{v=0}^{\infty} \underline{\mu}_{2v+u} Vec[\Omega_\omega] \tag{3.11}$$

$$where: \underline{\mu}_{2v} = E\{A^{\otimes 2v}(\omega)\} \ and \ \underline{\mu}_{2v+u} = E\{A(\omega)^v \otimes A(\omega)^{v+u}\} \tag{3.12}$$

**Proof.** Since the eigenvalues of the matrix $A(\omega)^{\otimes 2}$ have the form $\lambda_i \lambda_j$, $i,j = 1,...,p$ where $\lambda_i$ and $\lambda_j$ are the eigenvalues of the matrix $A(\omega)$ which its eigenvalues within the unit circle i.e. $|\lambda_i \lambda_j| < 1$ then $\left[I_{p^2} - A(\omega)^{2\otimes}\right]$ is not singular matrix,

$$
\begin{aligned}
Vec\left[\Upsilon_\omega(0)\right] &= Vec\left[E\left\{\underline{Y}_t(\omega)\underline{Y}'_t(\omega)\right\}\right] = E\left\{\left[I_{p^2} - A(\omega)^{\otimes 2}\right]^{-1} Vec\left[\Omega_\omega\right]\right\} \\
&= E\left\{\left[I_p^{\otimes 2} + \left[A(\omega)^{\otimes 2}\right] + \left[A(\omega)^{\otimes 2}\right]^2 + ...\right] Vec\left[\Omega_\omega\right]\right\} \\
&= \left[E\left\{I_p^{\otimes 2}\right\} + E\left\{A(\omega)^{\otimes 2}\right\} + E\left\{A(\omega)^{\otimes 4}\right\} + E\left\{A(\omega)^{\otimes 6}\right\}...\right] Vec\left[\Omega_\omega\right] \\
&= \sum_{v=0}^{\infty} \underline{\mu}_{2v} Vec\left[\Omega_\omega\right]
\end{aligned}
$$

therefore

$$
\begin{aligned}
vec\left[\Upsilon_\omega(u)\right] &= E\left\{\left(I_p \otimes A^u(\omega)\right)\left\{I_p^{\otimes 2} + A(\omega)^{\otimes 2} + A^{\otimes 4}(\omega) + ...\right\} Vec\left[\Omega_\omega\right]\right\} \\
&= E\left\{I_p \otimes A^u(\omega) + A(\omega) \otimes A^{1+u}(\omega) + A^2(\omega) \otimes A^{2+u}(\omega) + ...\right\} Vec\left[\Omega_\omega\right] \\
&= \sum_{v=0}^{\infty} \underline{\mu}_{2v+u} Vec\left[\Omega_\omega\right], \ u \geq 1
\end{aligned}
$$

∎

**Definition 3.1** *The spectral density of* (1.1) *model is given by*

$$S(\lambda) = \frac{1}{2\pi} \sum_{u=-\infty}^{\infty} \Upsilon_\omega(u) e^{-i\lambda u}, \ \lambda \in \mathbb{R} \tag{3.13}$$

**Theorem 3.2** *When* $S(\lambda)$ *exist then*

$$Vec\left[S(\lambda)\right] = \frac{1}{2\pi}\left[I_{p^2} + \sum_{u=1}^{\infty}\left(\underline{\mu}_u + \underline{\mu}'_u\right)\right] \sum_{v=0}^{\infty} \underline{\mu}_{2v} Vec\left[\Omega_\omega\right], \ pour \ \lambda = 0 \tag{3.14}$$

**Proof.**
$$S(\lambda) = \frac{1}{2\pi}\sum_{u=-\infty}^{\infty} \Upsilon_\omega(u) e^{-i\lambda u} = \frac{1}{2\pi}\left[\Upsilon_\omega(0) + \sum_{u=1}^{\infty} \Upsilon_\omega(u) e^{-i\lambda u} + \sum_{u=1}^{\infty} \Upsilon'_\omega(u) e^{i\lambda u}\right]$$

$$
\begin{aligned}
Vec[S(\lambda)] &= \frac{1}{2\pi}\left[Vec[\Upsilon_\omega(0)] + \sum_{u=1}^{\infty} Vec[\Upsilon_\omega(u)] e^{-i\lambda u} + \sum_{u=1}^{\infty} Vec[\Upsilon'_\omega(u)] e^{i\lambda u}\right] \\
&= \frac{1}{2\pi}\left(\sum_{v=0}^{\infty} \underline{\mu}_{2v} + \sum_{u=1}^{\infty} \underline{\mu}_u \sum_{v=0}^{\infty} \underline{\mu}_{2v} e^{-i\lambda u} + \sum_{u=1}^{\infty} \underline{\mu}'_u \sum_{v=0}^{\infty} \underline{\mu}_{2v} e^{i\lambda u}\right) Vec[\Omega_\omega]
\end{aligned}
$$

where $\underline{\mu}'_u = E\{A^u(\omega) \otimes I_p\}$. For $\lambda = 0$, we have

$$Vec[S(\lambda)] = \frac{1}{2\pi}\left[I_{p^2} + \sum_{u=1}^{\infty}\left(\underline{\mu}_u + \underline{\mu}'_u\right)\right] \sum_{v=0}^{\infty} \underline{\mu}_{2v} Vec[\Omega_\omega]$$

∎

**Theorem 3.3** *It is necessary that the matrix $\left[I_{p^2} - I_p \otimes A(\omega)\right]^2$ is invertible, for $S(\lambda)$ to exist and be continuous.*

**Proof.** If $\left[I_{p^2} - I_p \otimes A(\omega)\right]^2$ is invertible matrix, this is equivalent that the eigenvalues of the matrix $I_p \otimes A(\omega)$ are inside the circle unit i.e. those of the matrix $A(\omega)$ are also inside the circle unit i.e. $|\lambda_i| < 1$, $1 \leq i \leq p$, then

$$
\begin{aligned}
\left(\left[I_{p^2} - I_p \otimes A(\omega)\right]^2\right)^{-1} &= \left[I_{p^2} - \left(2I_p \otimes A(\omega) - [I_p \otimes A(\omega)]^2\right)\right]^{-1} \\
&= I_p^{\otimes 2} + 2\left[I_p \otimes A(\omega)\right] + 3\left[I_p \otimes A(\omega)^2\right] + 4\left[I_p \otimes A(\omega)^3\right] + ... \\
&= \sum_{u=0}^{\infty} (u+1)\left[I_p \otimes A(\omega)^u\right]
\end{aligned}
$$

thus

$$
E\left\{\left(\left[I_{p^2} - A^{\otimes 2}(\omega)\right]^2\right)^{-1}\right\} = \sum_{u=0}^{\infty} (u+1)\underline{\mu}_u
$$

since

$$
\begin{aligned}
\|Vec[S(\lambda)]\| &= \left\|\frac{1}{2\pi}\left(\sum_{v=0}^{\infty}\left\{I_{p^2} + \sum_{u=1}^{\infty}\left(\underline{\mu}_u + \underline{\mu}'_u\right)\right\}\underline{\mu}_{2v}\right)Vec[\Omega_\omega]\right\| \\
&\leq \left\|\sum_{u=0}^{\infty}(u+1)\underline{\mu}_u Vec[\Omega_\omega]\right\| < \infty
\end{aligned}
$$

∎

## 3.2 Estimation and asymptotic behavior

**Definition 3.2** *For $0 \leq u \leq T$ we define*

$$\widehat{\Upsilon}_N(u) = \frac{1}{(T-u+1)N}\sum_{t=1}^{T+u}\sum_{\omega=1}^{N}\underline{Y}_t(\omega)\underline{Y}'_{t-u}(\omega) \tag{3.15}$$

*which is an unbiased estimator of $\Upsilon_\omega(u)$.*

**Theorem 3.4** *The n.s. condition that*

$$\lim_{N\to\infty} Vec\left[\widehat{\Upsilon}_N(u)\right] = Vec[\Upsilon_\omega(u)] \quad p.s \tag{3.16}$$

*for all $u = 0, 1, ..., T$ is that the matrix $A(\omega)$ has all its eigenvalues within the unit circle.*

**Proof.** Consider, first $Vec\left[\widehat{\Upsilon}_N(0)\right]$ noting that

$$\sum_{t=0}^{T}\underline{Y}_t(\omega)\underline{Y}'_t(\omega), \quad \omega = 1, 2, ...$$

are *iid* random variables for $\omega = 1, ..., N$. So in order that $Vec\left[\widehat{\Upsilon}_N(0)\right] \to Vec[\Upsilon_\omega(0)]$ p.s, it is necessary that

$$\left\|Vec\left[E\left\{\sum_{t=0}^{T}\underline{Y}_t(\omega)\underline{Y}'_t(\omega)\right\}\right]\right\| < \infty$$

thus $\left\|Vec\left[E\left\{\underline{Y}_t(\omega)\underline{Y}'_t(\omega)\right\}\right]\right\| < \infty$, $t = 0, ..., T$, which is equivalent that the matrix $A(\omega)$ has all its eigenvalues within the unit circle.

To prove sufficiency, we note that for (4.20), it is sufficient for all $u$:

$$\left\|Vec\left[E\left\{\sum_{t=0}^{T-u}\underline{Y}_t(\omega)\underline{Y}'_{t-u}(\omega)\right\}\right]\right\| < \infty$$

which is true by inequality of Schwarz, if $\left\|Vec\left[E\left\{\underline{Y}_t(\omega)\underline{Y}'_t(\omega)\right\}\right]\right\| < \infty$, $t = 0, ..., T$, and this is verified if the matrix $A(\omega)$ has all its eigenvalues within the unit circle. ∎

**Definition 3.3** *We define the estimators $\widehat{\Omega}_N$ and $\hat{\underline{\mu}}_N(u)$ of $\Omega_\omega$ and $\underline{\mu}_u = E\{I_p \otimes A^u(\omega)\}$ respectively by*

$$Vec\left[\hat{\Upsilon}_N(u)\right] = \hat{\mu}_N(u) Vec\left[\hat{\Upsilon}_N(0)\right]^{-1}$$

$$Vec\left[\widehat{\Omega}_N\right] = \left[I_{p^2} - \hat{\underline{\mu}}_N(2)\right] Vec\left[\hat{\Upsilon}_N(0)\right]$$

**Theorem 3.5** *If the matrix $A(\omega)$ has all its eigenvalues within the unit circle then*

$$\lim_{N\to\infty} Vec\left[\widehat{\Omega}_N\right] = Vec\left[\Omega_\omega\right], \text{ and } \lim_{N\to\infty} Vec\left[\hat{\underline{\mu}}_N(u)\right] = Vec\left[\underline{\mu}_u\right] \text{ a.s, for } u = 1, ..., T-2 \quad (3.17)$$

*Moreover to that*

$$\sqrt{N}\left(Vec\left[\hat{\Upsilon}_N(0) - \Upsilon_\omega(0)\right]', ..., Vec\left[\hat{\Upsilon}_N(T) - \Upsilon_\omega(T)\right]'\right) \quad (3.18)$$

*converges as $N \to \infty$ to $(T+1)p^2$-dimensional normal variate wiyh zero means and covariance matrix blocks $\Sigma$, having $(p+1, q+1)$th matrix*

$$\lim_{N\to\infty} NCov\left(Vec\left[\hat{\Upsilon}_N(p)\right], Vec\left[\hat{\Upsilon}_N(q)\right]\right)$$

*it is necessary and sufficient that the matrix $A(\omega)$ has all its eigenvalues within the unit circle.*

**Proof.** Necessary condition : we note that

$$\sum_{t=0}^{T}\underline{Y}_t(\omega)\underline{Y}'_t(\omega), \quad \omega = 1, 2, ...N$$

are *iid* random variables.

In order that $\sqrt{N}\left[\hat{\Upsilon}_N(0) - \Upsilon_\omega(0)\right]$ is asymptotically normal, it is necessary and sufficient that :

$$\left\|E\left\{\sum_{t=0}^{T}\underline{Y}_t(\omega)\underline{Y}'_t(\omega)\right\}^2\right\| < \infty.$$

According to the central limit theorem for the *iid* random variables, which is true if and only if $\left\|E\left\{\underline{Y}_t(\omega)\underline{Y}'_t(\omega)\right\}^2\right\| < \infty$, but this inequality is equivalent that the matrix $A(\omega)$ has all its eigenvalues within the unit circle.

Sufficient condition : now, for $\sqrt{N}\left[\hat{\Upsilon}_N(u) - \Upsilon_\omega(u)\right]$ is also asymptotically normal for any $u \geq 1$, it is necessary and sufficient that

$$\left\|E\left(\sum_{t=0}^{T-u}\underline{Y}_t(\omega)\underline{Y}_{t-u}(\omega)\right)^2\right\| < \infty$$

which is true if the matrix $A(\omega)$ has all its eigenvalues within the unit circle, using Schwarz inequality.

Therefore the that the matrix $A(\omega)$ has all its eigenvalues within the unit circle is necessary and sufficient condition for any combination of the $\sqrt{N}\left[\hat{\Upsilon}_N(u) - \Upsilon_\omega(u)\right]$, $u \geq 0$, to be asymptoticallysoit multivariate normal. ∎

# 4   Bibliographie

# References

[1] J. Anděl. On Autoregressive Models With Random Parameters. Asymptotic Statistics (2). Edited by Peter Mandel and Marie Huškovà North-Holland (1984).

[2] C. Hsiao and M; Hashem Persaran.. Random Coefficient Panel Data Models. IZA DP N° 1236. August (2004).

[3] C. Hsiao. Analysis of Panel Data. Economic society monographs. N°. 34, $2^{nd}$ Edition, New York: Cambridge University Press, (2003).

[4] James D. Hamilton. Time-series analysis. Princeton University Press, Princeton New Jersey, (1994).

[5] Lon-Mu Liu. Random coefficient first-order autoregressive models. J. of Econometrics, Vol. 13, 305-325, (1980).

[6] F. Nicholls and Barry G. Quinn. Random Coefficient Autoregressive Models: An Introduction.Springer- Verlag New York Inc, (1982).

[7] P. M. Robinson. Statistical Inference for a Random Coefficient Autoregressive Model. Scand. J. Statist 5 Vol. 5, p163-168, (1978).

[8] P. A. V. B. Swamy. Statistical Inference in a Random Coefficient Regression Models. Springer, New York, (1971).

[9] P. A. V. B. Swamy. Efficient Inference in a Random Coefficient Regression Models. Econometrica, Vol. 38. N°. 2, 311-323, (1970).

[10] Vesna M. Ćojbašić. A Random Coefficient Autoregressive model (RCAR(1) Model). Univ. Beograd. Publ. Elektrotehn. Fak. Ser. Mat. 13, pages 45-50, (2004).